\documentclass[11pt]{amsart}
\usepackage{amsfonts,amssymb,amsbsy,amsmath,amscd,amstext,amsthm}

\usepackage{graphicx}
\newtheorem{theorem}{Theorem}
\newtheorem{thm}[theorem]{Theorem}

\newtheorem{proposition}{Theorem}
\newtheorem{lem}[proposition]{Lemma}

\newtheorem{cor}[theorem]{Corollary}

\newtheorem{conjecture}{Theorem}
\newtheorem{conj}[conjecture]{Conjecture}

\theoremstyle{remark}
\newtheorem{rem}{Remark}
\newtheorem{remark}[rem]{Remark}

\def\L{\mathcal L}

\def\Z{\mathbb Z}

\def\R{\mathbb R}

\begin{document}

\title{Discretized rotation has infinitely many periodic orbits}
\author{Shigeki Akiyama}
\address{Institute of Mathematics, University of Tsukuba, Tennodai 1-1-1,
Tsukuba, Ibaraki, 305-0006 Japan}
\email{akiyama@math.tsukuba.ac.jp}
\author{Attila Peth\H o}
\address{Department of Computer Science, University of Debrecen, P.O. Box 12, H-4010 Debrecen, Hungary}
\email{petho.attila@inf.unideb.hu}
\date{}
\thanks{The authors are supported by the Japanese Society for the Promotion of Science (JSPS), Grant in aid 21540010 and Invitation Fellowship Program FY2011, L-11514.}

\begin{abstract}
For a fixed $\lambda\in (-2,2)$, the discretized rotation on $\Z^2$ is defined by
$$
(x,y) \mapsto (y, -\lfloor x+\lambda y\rfloor).
$$
We prove that this dynamics has infinitely many periodic orbits.
\end{abstract}

\maketitle

\section{Introduction}
Space discretization of dynamical systems
attracted considerable interests among researchers \cite{Blank:94,Vladimirov:96,Vivaldi:94,
KKPV:97, BertheNouvel:07}.
One motivation is to understand the distance between original dynamics and
its computer simulation through discretized model.
In this paper, we are interested in a discretized planar rotation. It is very simple but we know surprisingly little on this discretized system.
We start with a conjecture studied by many authors, for e.g., in
\cite{Lowenstein-Hatjispyros-Vivaldi, Vivaldi:06, BLPV:03}
and from a point of view of shift radix system in \cite{Akiyama-Brunotte-Pethoe-Steiner:06}.
\begin{conj}
For all fixed $-2<\lambda<2$, all integer sequences $(a_n)$ defined by
\begin{equation}
\label{SRS2}
0\le a_{n+2}+\lambda a_{n+1}+a_n <1
\end{equation}
with initial value $(a_0,a_1)\in \Z^2$ are periodic.
\end{conj}
For $(x,y)=(a_n, a_{n+1})$, we have
$$
(a_{n+1},a_{n+2})=(y,-\lfloor x+\lambda y\rfloor),
$$
and it defines a map $F:(x,y)\mapsto (y,-\lfloor x+\lambda y\rfloor)$ on
$\Z^2$.
In other words, we are interested in the dynamics $F$ on $\Z^2$:
\begin{equation}
\label{Algo}
\begin{pmatrix} x\\ y\end{pmatrix} \mapsto
\begin{pmatrix} 0&1 \\
                -1 & -\lambda
\end{pmatrix}
\begin{pmatrix} x\\y\end{pmatrix}+
\begin{pmatrix} 0\\ \langle \lambda y \rangle \end{pmatrix},
\end{equation}
where $\langle x \rangle= x-\lfloor x\rfloor$.
Since the eigenvalues of the matrix are two conjugate complex numbers of modulus one,
this dynamics can be regarded as a rotation
having invariant confocal ellipses, acting on the lattice $\Z^2$.
After the rotation of
angle $\theta$ with $\lambda=-2\cos\theta$,
we translate by a small vector to make the image lie in $\Z^2$.
An affine equivalent formulation using Euclidean rotation is found in the next section.
From the shape of the inequality in (\ref{SRS2}), the dynamics (\ref{Algo})
is reversible, i.e., $F(x,y)=(y,z)$ implies $F(z,y)=(y,x)$.
Thus we have $\phi F^{-1}= F \phi$ with $\phi(x,y)=(y,x)$ and
$F$ is a bijection on $\Z^2$. In other words, $F$ is a composition of
two involutions $F \phi$ and $\phi$. 

The conjecture is supported by
numerical experiments \cite{Vivaldi:94, Akiyama-Brunotte-Pethoe-Steiner:06}. It is also expected from a heuristic ground:
cumulation of errors of $F^n$ from
the exact $n\theta$ rotation is expected to be small and seemingly impossible to
avoid hitting the same lattice points.
The cumulative error bound is discussed in \cite{KKPV:97,Vivaldi-Vlasimirow:03}.
However this problem is notorious, and our knowledge
is limited. We only know that Conjecture 1 holds for
11 values $\lambda=0,\pm 1,(\pm 1 \pm\sqrt{5})/2, \pm \sqrt{2},\pm \sqrt{3}$,
see \cite{Akiyama-Brunotte-Pethoe-Steiner:07, Lowenstein-Hatjispyros-Vivaldi, Akiyama-Brunotte-Pethoe-Steiner:06}.
Apart from three trivial cases $0, \pm 1$, the proof is
highly non trivial and uses the self-inducing structure found
in the associated planer piecewise isometry when $\theta/\pi$ is rational and
$\lambda$ is quadratic.
If $\theta/\pi$ is rational, then we can embed the problem into
piecewise isometry acting on a certain higher dimensional torus
(see \cite{Lowenstein-Hatjispyros-Vivaldi, Kouptsov-Lowenstein-Vivaldi:02}, and
also \cite{Ashwin:97, Ashwin-Chambers-Petrov:97} for connection to digital filters).
Piecewise isometries have zero entropy \cite{Buzzi:01},
but we know little on their periodic orbits \cite{Goetz:02}.
It is noteworthy that a certain piecewise isometry generated by $7$-fold rotation in
the plane is governed by several self-inducing structures
\cite{Lowensein-Koupstov-Vivaldi:04, GoetzPoggiasparlla:04, Akiyama-Harriss:12},
but it is irrelevant to the map $F$.
If $\lambda$ is a
rational number whose denominator is the power of a prime, then the dynamics
is understood as the composition of $p$-adic rotation and symbolic shift
in \cite{Bosio-Vivaldi:00},
but it seems difficult to extract information
on periodic orbits through this embedding.
At this stage, we are interested in
giving a non trivial general statement for this dynamics.
In this note, we will show

\begin{thm}
\label{Main}
For all fixed $\lambda\in (-2,2)$ there are infinitely many periodic orbits
of the dynamics (\ref{Algo}) on $\Z^2$.
\end{thm}

\noindent
More precisely,
we prove that there are infinitely many {\it symmetric} periodic orbits
(see \S 3 for the definition). 
See Corollary \ref{Low} in \S \ref{Gen},
for a qualitative statement.
Theorem \ref{Main} is
new for all $\lambda$ except the above $11$ values, and
gives another support of the conjecture.
Note that the idea of observing symmetric periodic orbits
dates back to G.~Birkhoff, who showed the existence of infinitely many symmetric
periodic orbits for the restricted three body problem \cite{Birkhoff:15, Birkhoff:66},
whose dynamics is composed of two involutions.

Adding a counting technique of lattice points in number theory,
we can generalize Theorem \ref{Main} to the sequences generated by:
$$
-\eta \le a_{n+2}+\lambda a_{n+1}+a_n <1-\eta
$$
for a fixed $\eta\in \R$, as in \S \ref{Gen}, Theorem \ref{Main2}.
It covers a significant class of discretized rotations, but
we do not know how large is this two parameter family
within the set of all
invertible discretized rotations in $\Z^2$, up to conjugacy.

We say $p=p(x,y)>0$ is the {\it period} of $(F^n(x,y))_{n\in \Z}$, if
it is the smallest positive integer $p$ with $F^{p}(x,y)=(x,y)$.
If there is no such $p$, then $p(x,y)$ is not defined.
It is remarkable that
the distribution of periodic orbits drastically changes by whether
$\theta/\pi$ is irrational or rational. We have

\begin{lem}
\label{Period}
Let $\theta/\pi$ be irrational and $p$ is a positive integer. Then
 there are only finitely many periodic orbits of period $p$.
\end{lem}

This fact follows from Theorem 2.1 of \cite{Vivaldi:94} but we give a quick proof
in \S 5.
On the other hand, if $\theta/\pi$ is rational,
in view of the above torus embedding,
it is natural to obtain infinitely
many periodic orbits of period $p$, which fall into the same period cell.
Theorem 4.3 in \cite{Akiyama-Brunotte-Pethoe-Steiner:06} gives
a concrete example of infinite periodic orbits of period $p$
for $\theta=(1-1/p)\pi$ with an odd prime $p$.

By Theorem \ref{Main} and Lemma \ref{Period}, we know that
there exist arbitrary large periods, if $\theta/\pi$ is irrational.
We expect the same to hold for all $\lambda\neq 0,\pm 1$, but there are proofs only
for the above eight quadratic cases.

\section{Setting and strategy}

\begin{figure}[h]
\includegraphics[width=7cm]{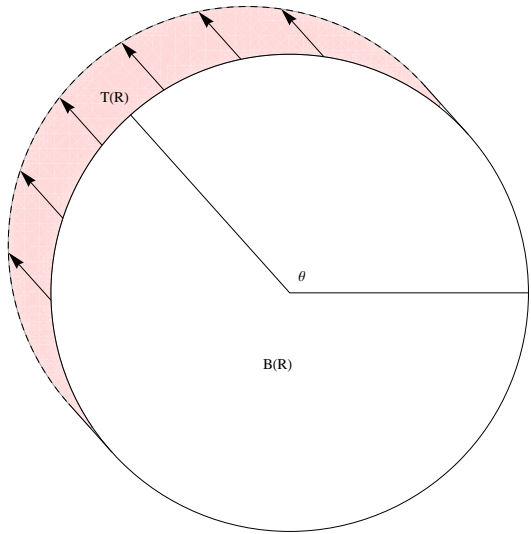}
\caption{Trap Region\label{Trap}}
\end{figure}

Let
$\lambda=-2\cos\theta$ where $\theta$ is a real number in $(0,\pi)$
and let
$Q=\begin{pmatrix}-\sin \theta& \cos \theta\\ 0 & 1\end{pmatrix}$.
Our transformation on $\Z^2:(x,y)\mapsto (X,Y)$
is written as
$$
\begin{pmatrix} X\\ Y\end{pmatrix} =
\begin{pmatrix} 0&1 \\
                -1 & -\lambda
\end{pmatrix}
\begin{pmatrix} x\\y\end{pmatrix}+
\begin{pmatrix} 0\\ \mu \end{pmatrix}
$$
with $\mu\in [0,1)$.
Since
$$
Q \begin{pmatrix} \cos \theta & -\sin \theta\\
                  \sin \theta & \cos \theta\end{pmatrix} Q^{-1}
=
\begin{pmatrix} 0 & 1 \\-1 & -\lambda
\end{pmatrix},
$$
we view this algorithm as
\begin{equation}
\label{Algo2}
Q^{-1} \begin{pmatrix} X\\ Y\end{pmatrix} =
\begin{pmatrix} \cos \theta & -\sin \theta \\
                \sin \theta & \cos \theta
\end{pmatrix} Q^{-1}
\begin{pmatrix} x\\y\end{pmatrix}+ Q^{-1}
\begin{pmatrix} 0\\ \mu \end{pmatrix}.
\end{equation}
Thus it is the dynamics acting on the lattice
$\L=\begin{pmatrix} -\csc \theta \\0 \end{pmatrix} \Z+ \begin{pmatrix}
\cot \theta\\ 1\end{pmatrix} \Z$
written as
the composition of the Euclidean rotation of angle $\theta$
followed by a small translation
$$
\mathbf{v} \mapsto \mathbf{v}+\mu \begin{pmatrix} \cot \theta\\1\end{pmatrix}
$$
with $\mu \in [0,1)$.
Let $R$ be a positive real number and $B(R)$ be a ball of radius $R$ centered at
the origin. Define a {\it trap region} $T(R)$ by
$$
T(R)= \left\{ \left. x+y \begin{pmatrix}
\cot \theta\\ 1\end{pmatrix} \ \right| x\in B(R), y\in [0,1)\right\}
\setminus B(R).
$$
The situation is demonstrated in Figure \ref{Trap}.

Now we explain the strategy of the proof.
It is clear from the description of the dynamics,
that every unbounded orbit
starting from a point in $\L\cap B(R)$ must visit at least
once the trap region $T(R)$.
Assume that there are only finitely many periodic orbits of (\ref{Algo2}).
Since we are dealing with dynamics on the lattice $\L$,
periodicity of an orbit is equivalent to its boundedness.
We
argue by contradiction: the assumption that there are just a finite
number of periodic orbits would imply that, for $R \gg 1$, the number
of lattice points in $T(R)$ are strictly less than the number of symmetric
unbounded orbits starting from $B(R)$, which is impossible.

\section{Lower bound of unbounded orbits}

Symmetric
periodic orbits of time-reversal dynamics had been studied
for a long time.
We shall make
use of a well-known property which may belong to a folklore.
Let $(a_n)$ be a bi-infinite integer sequence and $b$ be an integer.
We say that $(a_n)$ is periodic, if
$a_{n+b}=a_n$ hold for all $n$. The sequence $(a_n)$
is {\it symmetric} at $b/2$,
if $a_{b-n}=a_n$ holds for all $n$, and it is {\it doubly symmetric} if
it is symmetric at $b_1/2$ and $b_2/2$ with $b_1\neq b_2$.
%We refer to such a point, like $b/2$, the {\it center} of symmetry.

\begin{lem}
\label{Sym}
If a sequence is doubly symmetric then
it is periodic. Moreover symmetric and periodic sequences are doubly symmetric.
\end{lem}

\proof Assume that $(a_n)$ is symmetric at $b_1$ and $b_2$ with $b_1\neq b_2$.
Then $a_{n+b_2-b_1}=a_{b_1-n}=a_n$. Let $(a_n)$ be periodic of period $b$ and
symmetric at $c$. Then $a_{c-b-n}=a_{n+b}=a_n$.
\qed
\bigskip

In fact, this lemma is a restatement of
Theorems 1 and 2 in \cite{DeVogelaere}, where
symmetric periodic orbits of time
reversal dynamics composed of two involutions are studied.
The reader finds precise description on such periodic orbits there.

Let $(x,y)\in \Z^2$. To the bi-infinite orbit $(F^n(x,y))_{n\in \Z}$ we can associate
uniquely the bi-infinite sequence $(a_n)$ consisting of the 1-st coordinates of the elements of the orbit.
It is clear that
$(F^n(x,y))$ is periodic if and only if $(a_n)$ is periodic.
%We also say that an orbit $(F^n(x,y))$ is symmetric if $(a_n)$ is so.
%We note that symmetric
%periodic orbits of time-reversal dynamics has been studied already
%in old literatures, see \cite{DeVogelaere}.
Hereafter we identify the orbit $(F^n(x,y))$ and the bi-infinite sequence
$(a_n)$ and say that an orbit $(F^n(x,y))$ is symmetric if $(a_n)$ is so.
Assume that $(a_n)$ is symmetric at $b/2$.
If $b$ is odd,
then $a_{(b-1)/2}=a_{(b+1)/2}$ and the orbit is of the form:
$$
\dots,c_3,c_2,c_1,X,X,c_1,c_2,c_3,\dots
$$
with some $X\in \Z$ and a sequence $(c_n)\subset \Z$. Clearly $(c_n)$ is determined by $X$.
Such orbits are in $\mathrm{Fix}(\phi)$, the set of orbits fixed by the involution
$\phi$.
If $b$ is even, then $a_{b/2-1}=a_{b/2+1}$ and the orbit is of the form
$$
\dots,c_3,c_2,c_1,X,Y,X,c_1,c_2,c_3,\dots
$$
for some $X, Y\in \Z$ and $(c_n)\subset \Z$. Of course
$(c_n)$ is determined by $X$ and $Y$.
These orbits belong to $\mathrm{Fix}(F\phi)$, the set fixed by the other involution
$F\phi$.
By Lemma \ref{Sym}, each periodic orbit of $(F^n(x,y))$
is doubly symmetric.
It belongs to the same fix point set
if the period is even, and to the different fix point sets if
the period is odd.

\begin{remark}
Not all orbits are symmetric.
For e.g., if $\lambda=(1+\sqrt{5})/2$ then we have
$$
(-1,4)\rightarrow
(4,-6)\rightarrow
(-6,5)\rightarrow
(5,-3)\rightarrow
(-3,-1)\rightarrow
(-1,4).
$$
We do not know a way to estimate from below the number of
asymmetric orbits.
\end{remark}

By our
assumption of
reductio ad absurdum,
there exist only a finite number, say $C_1$, of periodic orbits.
For any $R$, the number of points in $\L\cap B(R)$
whose orbits are periodic is less than $C_1$.

\subsection{Unbounded orbits in $\mathrm{Fix}(\phi)$}

Let $(a_n)$ be an unbounded orbit in $\mathrm{Fix}(\phi)$.
By Lemma \ref{Sym}, $(a_n)$ can not be doubly symmetric, which 
implies that two orbits starting from different fixed points in $\mathcal{L}$
never intersect.
Apart from a finite number of exceptions,
points of the form
$$
X \begin{pmatrix} -\csc \theta \\0 \end{pmatrix} +
X \begin{pmatrix} \cot \theta\\ 1\end{pmatrix} \in \L\cap B(R)
$$
generate distinct unbounded orbits.
From
$$
X^2 \left(-\csc \theta+\cot \theta \right)^2 + X^2 \leq R^2,
$$
we conclude that there are at least $2R \cos(\theta/2)-C_1$ unbounded
orbits in $\mathrm{Fix}(\phi)$ starting from $\L\cap B(R)$.

\subsection{Unbounded orbits in $\mathrm{Fix}(F \phi)$}

Similarly, by Lemma \ref{Sym}, $(a_n)$ can not be doubly symmetric.
So, by (\ref{SRS2}),
our task is to count the number of the pairs $(X,Y)$ which satisfy
\begin{equation}
\label{Band}
0\le X + \lambda Y + X < 1
\end{equation}
and
\begin{equation}
\label{Ball}
X \begin{pmatrix} -\csc \theta \\0 \end{pmatrix} + Y \begin{pmatrix}
\cot \theta\\ 1\end{pmatrix} \in B(R).
\end{equation}
For this computation, we substitute the inequality (\ref{Band}) by
\begin{equation}
\label{Band2}
-1\le X + \lambda Y + X < 1
\end{equation}
and count the number of pairs $(X,Y)$ satisfying
(\ref{Band2}) and (\ref{Ball}).
It is clear that for a fixed $Y$, there is a unique $X$ which
satisfies (\ref{Band2}). Since $X=Y \cos \theta+\varepsilon$
with $|\varepsilon|\le 1/2$, we have
\begin{equation}
\label{YRange}
\left(Y \cot \theta - X \csc \theta \right)^2 + Y^2
=Y^2+\frac{\varepsilon^2}{\sin^2 \theta} \le R^2
\end{equation}
from (\ref{Ball}), and we have at least $2R-C_2$ such points with
a non negative constant $C_2$.

If $\lambda$ is irrational,
then there is no $(X,Y)\neq (0,0)$ which satisfies either
$$
-1= 2X + \lambda Y \quad \text{ or } \quad 0=2X+ \lambda Y.
$$
Using the symmetry $(X,Y)\leftrightarrow (-X,-Y)$,
we see that the number of $(X,Y)$ with (\ref{Ball}) and
$$
-1\le X + \lambda Y + X < 0
$$
is exactly
one less than the number of $(X,Y)$ with (\ref{Band}) and (\ref{Ball}), which
counts the origin.
Thus the number of $(X,Y)$ having (\ref{Band}) and (\ref{Ball})
is at least $R-C_2/2$.

If $\lambda$ is rational, then we additionally have to
take care of the points $(X,Y)$ on the line $-1= 2X + \lambda Y$
and  $0= 2X + \lambda Y$.
However we can easily show that
the number of $(X,Y)$
with (\ref{Ball}) on the line $-1= 2X + \lambda Y$ and the one
on the line $0= 2X + \lambda Y$ differ only by some constant.
Thus in any case,
there are at least $R-C_3$ unbounded
orbits in $\mathrm{Fix}(F\phi)$ starting from $\L\cap B(R)$
with a non negative constant $C_3$.

\section{Lattice points in the trap region}

By construction, if the trap region $T(R)$ and the line $x=y \cot \theta +c$
has non empty intersection, then it is a half-open interval of length $\csc \theta$.
Thus if the line $x=y \cot \theta +c$ intersects
$\L \cap T(R)$ then it is a single point.
The lattice $\L$ is covered by a family of parallel lines:
$$
\Xi=\{ x=y \cot \theta - k \csc \theta \ |\ k\in \Z \}.
$$
%A relation:
%$$
%\begin{pmatrix} -\csc \theta \\0 \end{pmatrix}
%=
%\begin{pmatrix} - \sin \theta\\ \cos \theta \end{pmatrix}
%-\cos \theta \begin{pmatrix} \cot \theta\\ 1\end{pmatrix},
%$$
%gives a decomposition of the vector $\begin{pmatrix} -\csc \theta \\0 \end{pmatrix}
%$
%into $\begin{pmatrix} \cot \theta\\ 1\end{pmatrix}$ direction
%and its perpendicular direction.
We easily see that the distance between adjacent lines of $\Xi$
is $1$. Thus there are exactly $2\lfloor R \rfloor+1$
points in $\L \cap T(R)$.
\begin{remark}
If $\theta<2\pi/3$, then
$2R\cos (\theta/2)-C_1+R-C_3 > 2R+1$ holds for sufficiently large
$R$ and we immediately obtain the desired contradiction.
\end{remark}

Let us take into account the reversibility of $F$. Since we are dealing
with unbounded symmetric orbits starting from $\L \cap B(R)$,
if $(a_n,a_{n+1})=(C,D)$ then there is
an index $m$ such that $(a_m,a_{m+1})=(D,C)$.
%, and $C\neq D$ except when
%it is symmetric of $(X,X)$ type at $n+1/2$.
Let $\Phi: \; \mathcal{L} \mapsto \mathcal{L}$ be defined as follows
$$
\Phi:  x \begin{pmatrix} -\csc \theta \\0 \end{pmatrix} +
y \begin{pmatrix} \cot \theta\\ 1\end{pmatrix} \mapsto
 y \begin{pmatrix} -\csc \theta \\0 \end{pmatrix} +
x \begin{pmatrix} \cot \theta\\ 1\end{pmatrix}.
$$
If an orbit visits $\Phi(T(R))\cap T(R)$ then the number of visits is at least two.
In other words, we only have to count the number of lattice points up to
this symmetry by $\Phi$ in $T(R)$.

The mapping $\Phi$ is the reflection with respect to
the vector $\begin{pmatrix} -\csc \theta+\cot \theta \\ 1 \end{pmatrix}
$,
because the two vectors
$\begin{pmatrix} -\csc \theta \\0 \end{pmatrix}$
and $\begin{pmatrix} \cot \theta\\ 1\end{pmatrix}$
have the same length.
Thus the reflection $\Phi$ leaves the vector
$\begin{pmatrix} -\sin (\theta/2) \\ \cos (\theta/2) \end{pmatrix} = \begin{pmatrix} -\csc \theta+\cot \theta \\ 1 \end{pmatrix} \cos (\theta/2)$
invariant.
% which is perpendicular to
%$\begin{pmatrix} \cos (\theta/2) \\ \sin (\theta/2) \end{pmatrix}$.

\begin{figure}[h]
\includegraphics[width=8cm]{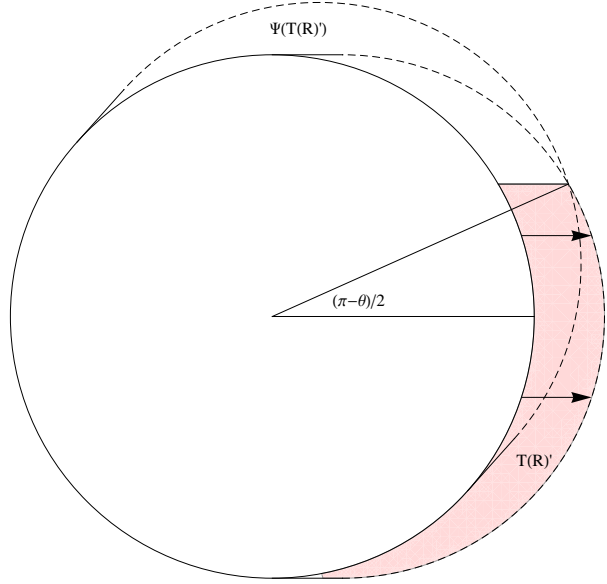}
\caption{Symmetry of the trap region\label{Count}}
\end{figure}

To make computation easy, we rotate $T(R)$ and $\L$ by $-\theta$ and present
the situation in Figure \ref{Count}.
$T(R)'$, $\L'$ are the images by
this rotation and $\Psi$ is the corresponding reflection. Then every line
$y=k$ with $k\in \Z \cap [-R,R]$ contains a single point in $\L'\cap T(R)'$ and the reflection $\Psi$ leaves the vector
$$
\begin{pmatrix} \cos \theta & \sin \theta \\
                -\sin \theta & \cos \theta
\end{pmatrix} \begin{pmatrix} -\sin(\theta/2) \\ \cos (\theta/2) \end{pmatrix} = \begin{pmatrix} \sin (\theta/2) \\ \cos (\theta/2) \end{pmatrix}
$$
invariant.

So our task is to estimate from above the
number of lattice points in $\L'\cap T(R)'$ which are below the line
$y= x \cot (\theta/2)$.
The intersection of the line and the boundary of $T(R)'$ with
the largest $y$-coordinate is
$$
\left(\sqrt{R^2 \sin^2 \left(\frac {\theta}2\right)-\frac 14}+\frac12 \tan
   \frac{\theta}{2},\cot \frac{\theta}{2}
   \sqrt{R^2 \sin^2 \left(\frac {\theta}2\right)-\frac 14}+\frac 12\right)
$$
and we have
$$
\cot \frac{\theta}{2}
   \sqrt{R^2 \sin^2 \left(\frac {\theta}2\right)-\frac 14}+\frac 12
   =R \cos \left(\frac{\theta}{2}\right)+\frac{1}{2}+ O\left(\frac 1R\right).
$$
We count the number of points whose $y$-coordinate does not exceed this value,
i.e., the points in the shaded part in Figure \ref{Count}.
Thus the number of lattice points up to symmetry
in $\L\cap T(R)$ is bounded from above by $R+R\cos (\theta/2)+C_4$ with
a non negative constant $C_4$.

\section{Proof of Theorem \ref{Main} and Lemma \ref{Period}}

From the assumption that there are only finitely many periodic orbits,
we derived several estimates in the previous sections.
By Lemma \ref{Sym},
unbounded orbits in $\mathrm{Fix}(\phi)$ and those in $\mathrm{Fix}(F\phi)$ have no
intersection. Thus
$$
2R\cos (\theta/2) -C_1 + R-C_3
$$
distinct unbounded orbits must visit $T(R)$ and there are
only
$$
R+R\cos (\theta/2)+C_4
$$
lattice points in $\L\cap T(R)$ up to symmetry. However
$$
2R\cos (\theta/2) -C_1 + R-C_3 \le R+R\cos (\theta/2)+C_4
$$
does not hold for sufficiently large $R$. The proof of Theorem \ref{Main}
is finished.

Let us show Lemma \ref{Period}. 
%First
%consider the case $p=2$. The periodic orbit of period 2
%is of the form:
%$$
%(x,y)\rightarrow (y,x)\rightarrow (x,y)
%$$
%and it is easy to see that there are only finitely many $(x,y)$ which satisfy
%$$
%0\le x+\lambda y+x<1 \quad  \text{ and } \quad
%0\le y+\lambda x+y<1.
%$$
Assume that there are infinitely many $(x,y)$ that
$F^{p}(x,y)=(x,y)$ with $p>2$. By induction using (\ref{Algo2}), we have
$$
\begin{pmatrix} u\\v\end{pmatrix}
= \begin{pmatrix} \cos p \theta& -\sin p \theta\\
                          \sin p \theta & \cos p \theta\end{pmatrix} \begin{pmatrix} u\\v\end{pmatrix}+ \sum_{i=1}^{p} \mathbf{v}_i
$$
where $\begin{pmatrix} u\\v\end{pmatrix}\in \L$ and
$\Vert \mathbf{v}_i \Vert \le \csc \theta$. Here $\Vert \cdot \Vert$ is the
Euclidean norm.
However if $\begin{pmatrix} u\\v\end{pmatrix}$ is sufficiently large, then
$$
\left\Vert \begin{pmatrix} u\\v\end{pmatrix}
- \begin{pmatrix} \cos p \theta& -\sin p \theta\\
                          \sin p \theta & \cos p \theta\end{pmatrix} \begin{pmatrix} u\\v\end{pmatrix}
\right\Vert\ge c
\left\Vert 
                          \begin{pmatrix} u\\v\end{pmatrix}
\right\Vert>
p \csc \theta.$$ 
Here $c$ is computed as the operator norm:
$$
c=\frac{1}{\left\Vert 
\begin{pmatrix} 1-\cos p \theta& \sin p \theta\\
                          -\sin p \theta & 1-\cos p \theta\end{pmatrix}^{-1}             
\right\Vert}=2\left|\sin\left(\frac{p \theta}2\right)\right|
$$
which is positive since $\theta/\pi$ is irrational. 
%Here
%we use the fact $p>2$ and replace $(x,y)$ with $F(x,y)$ to make the left side large,
%if it is necessary. 
This gives a contradiction.

\section{Generalization}
\label{Gen}
One can generalize the result to the sequences defined by:
$$
-\eta\le a_{n+2}+\lambda a_{n+1}+a_n <1-\eta,
$$
with $\eta\in \R$.
This kind of interval shifts are studied, for e.g., in
\cite{Kouptsov-Lowenstein-Vivaldi:02, Surer:09}.
In complete analogy to our main result, we have
\begin{thm}
\label{Main2}
For a fixed $\lambda\in (-2,2)$ and $\eta\in \R$,
there are infinitely many periodic orbits
of the dynamics 
\begin{equation}
\label{SRSeta}
(x,y)\rightarrow (y, -\lfloor \lambda y+x+\eta \rfloor)
\end{equation}
on $\Z^2$.
\end{thm}

We can deduce a qualitative statement:

\begin{cor}
\label{Low}
There is a positive constant $C$ depending on $\lambda$ and $\eta$ such
that within $B(R)$, there are at least $CR$ periodic orbits
of (\ref{SRSeta}).
\end{cor}

Hereafter we sketch the
proof of Theorem \ref{Main2} and Corollary \ref{Low}.
Putting $\kappa=\eta/(2+\lambda)$, the inequality becomes
$$
0\le (a_{n+2}+\kappa)+\lambda (a_{n+1}+\kappa)+(a_n+\kappa) <1.
$$
Therefore by substituting $\L$ with
$\L'=\L+Q^{-1}\begin{pmatrix}\kappa\\ \kappa\end{pmatrix}$,
our algorithm has exactly the same shape as (\ref{Algo2}).
Though the error term becomes worse than the one in \S 4,
we can
show that the number of lattice points of $\L'$ up to symmetry
within the trap region is
$$
R+ R\cos (\theta/2)+O(R^{2/3+\epsilon})
$$
for any positive constant $\epsilon$.
Here we used the method of
Vinogradov to count the number of lattice points
in the cylindrical region bounded by curves of positive curvature,
for e.g., see p.8-22 of \cite{Vinogradov:85} or
\cite{Vinogradov:54, Huxley:87}.

Similarly to \S 3, there are $2R \cos (\theta/2)-C_1$ unbounded orbits in
$\mathrm{Fix}(\phi)$.
We count unbounded orbits in $\mathrm{Fix}(F\phi)$, i.e., the
number of $(X,Y)$, which satisfy:
\begin{equation}
\label{UD}
\lambda Y/2 \bmod{1} \cap [-\eta/2, (1-\eta)/2) \neq \emptyset
\end{equation}
and (\ref{YRange}).

If $\lambda$ is irrational, then $(\lambda/2)Y \bmod{1}$
is uniformly distributed and the number of such $Y$'s is $R+o(R)$.
When $\lambda$ is rational, put
$\lambda/2=p/q$ with $(p,q)=1$. Then $(\lambda/2)Y \equiv i/q \bmod{1}$ for $i\in \{0,1,\dots, q-1\}$
with the same frequency $1/q$.

Let us study the case that $q$ is even.
Since
$$
\{i/q \bmod{1}\} \cap [-\eta/2, (1-\eta)/2)
$$
has cardinality $q/2$,
the number of points with (\ref{UD}) and (\ref{YRange})
is again $R+o(R)$. Once we have this estimate $R+o(R)$ then
$$
2R \cos(\theta/2)-C_1+R+o(R)\le R +R \cos(\theta/2)+o(R^{2/3+\epsilon})
$$
does not hold for sufficiently large $R$. Here $C_1$ appears only once
in the left side, because it is the number of periodic orbits of the system.
We obtain the contradiction for $R\gg 1$.

It remains to show the case when $q$ is odd. Then
$$
\{i/q \bmod{1}\} \cap [-\eta/2, (1-\eta)/2)$$
has cardinality either $(q-1)/2$ or $(q+1)/2$ depending on $\eta$.
Thus the number of unbounded orbits in $\mathrm{Fix}(F\phi)$
is bounded from below
by $R-R/q+o(R)$.
Thus we have to show that
$$
2R \cos(\theta/2)-C_1+R-R/q+o(R)>R +R \cos(\theta/2)+o(R^{2/3+\epsilon})
$$
for large $R$. This is valid because
$$
\cos(\theta/2) > 1/q
$$
holds for $q>2$, since $\cos (\theta/2)=\sqrt{(1+\cos(\theta))/2}=
\sqrt{(1-p/q)/2}$.

The validity of statement of Corollary \ref{Low} is
invariant under affine transformations fixing the origin, up to
appropriate changes the constant $C$.
To prove Corollary \ref{Low},
note that we may take
$$
C_1= (1-\varepsilon) R \left(\cos (\theta/2)-1/q\right)
$$ with
a small $\varepsilon>0$ in the above proof
to get the same contradiction.
Here the term $-1/q$ is necessary only in the last case of the proof.

\section{Acknowledgments}
We are grateful to Franco Vivaldi for giving us comments and references
of the earlier version of this article.
The first author moved from Niigata University to the current address in August 2012.
The paper was written, when the second author was visiting Niigata University as
a long term research fellow of JSPS. Both of us wish to express our deep gratitude
to all the staffs in Department of Mathematics, Niigata University and
the support from JSPS.

%\bibliographystyle{amsplain}
%\bibliography{reflist}

\end{document}